\def\tpt{\hskip 20pt}

\input amstex
\input amsppt.sty
\magnification=\magstep1
\hsize=15truecm
\vsize=20truecm
\TagsOnRight
\pageno=1

{\tpt}

\vskip 20pt

\centerline {\bf A Direct Formulation of Dieudonn\'e Module
Theory\footnote""{Mathematics Subject Classification
(1991): 14H05, 14H10, 14G27}}

\vskip 10pt \centerline {Kezheng Li\footnote"*"{Supported by NSFC
of China, grant number 11031004}}

\vskip 10pt
\centerline{\sl Department of Mathematics, Capital
Normal University, Beijing 100048, China}

\vskip 20pt
\flushpar{\bf Abstract.} We define a Dieudonn\'e module as the module of
Dieudonn\'e elements, and set up Dieudonn\'e module theory in a simple way.
Under this formulation we give explicit formulae for the duality and the corresponding
differential operators.

\vskip 20pt
\flushpar {\bf 0. Introduction}

\vskip 10pt
The Dieudonn\'e module theory of commutative group schemes has been set up for
more than fourty years (see [Dem] and [Laz]). It gives an embedding of the category of
finite commutative
group schemes over a perfect field $k$ to the category of modules over the ring
$$A=W[F,V]/(FV-p,VF-p,Fa-a^\sigma F,Va-a^{\sigma^{-1}}V\ \forall a\in W)$$
(where $W=W(k)$, the ring of Witt vectors over $k$) as a full
subcategory. Since then the study of commutative group schemes, especially over
a field of characteristic $p>0$, have almost been the study of Dieudonn\'e modules.

However, for a long time the ways to formulate Dieudonn\'e module theory had been quite indirect.
Basically there were two ways: One was to set up Dieudonn\'e module theory for
$p$-divisible groups first, and then show that every finite commutative group scheme
can be embedded into a $p$-divisible group as a closed subgroup scheme (see [Dem], [dJ1]
or [Laz]); the other used crystalline cohomology (see [BBM] or [Me]). Of course these
ways of formulation are by no means simple.

On the other hand, there was a simple definition of the Dieudonn\'e module functor $D$,
namely for a finite commutative group scheme over $k$
such that both $G$ and $G^D$ are
infinitesimal, define $D(G)=\text{Hom}_k(G,\Cal W)$, where $\Cal W$ is the Witt scheme over $k$
(see [Dem] or [LO]).
However, the proofs of the faithfulness of $D$ were still quite indirect and by no means simple
(see e.g. [dJ1]). Until recently we saw a direct proof of this by Pink (see [P]).

In this paper we follow this idea, and we define $D(G)$ to be a subset of ``Dieudonn\'e
elements'' in the structure ring of $G$. This is a simple and direct formulation.

For preparations, in \S 1 we write down some basic facts about differential operators,
and in \S 2 we write down some basic facts about the Witt scheme.
For the convenience of the reader, we write down complete proofs of the well-known
facts. Then in \S 3 we give a direct formulation of Dieudonn\'e module
theory. In particular we give a simple proof of the faithfulness of $D$.

Under this formulation, it is natural to ask what are the Dieudonn\'e
elements of the Cartier dual. This problem is equivalent to giving an explicit isomorphism
of the truncated Witt scheme $\Cal W_{n,n}$ to its Cartier dual (under the usual
coordinates). We give an answer in \S 4 for this problem: under the usual isomorphism
$$\Cal W_{n,m}\cong\text{Spec}(A_{n,m})$$
where
$$A_{n,m}=\Bbb F_p[x_0,...,x_n]/(x_0^{p^{m+1}},...,x_n^{p^{m+1}})$$
(and $x_0,...,x_n$ are {\bf specially chosen} in the definition of Witt scheme),
the following $\Bbb F_p$-linear functional on $A_{n,m}$
$$y(\prod_{i=0}^nx_i^{j_i})=\cases 1 & \text{ if } j_1=...=j_n=0,\ j_0=p^m \cr
0 & \text{ otherwise} \endcases$$
is a ``standard'' Dieudonn\'e element of $\Cal W_{n,m}^D\cong\text{Spec}(A_{n,m}^D)$, which gives
a standard isomorphism $\Cal W_{n,m}^D\to\Cal W_{m,n}$. This fact is by no means obvious.

In \S 5 we give an explanation of Dieudonn\'e elements as differential operators. 
On $\Bbb Z[x]$ we can define differential operators
$$D^{(p^i)}=\frac{1}{p^i!}\frac{\text{\rm d}^{p^i}}{\text{\rm d}x^{p^i}} \tpt (0\le i\le r)$$
which can also be defined on $R=\Bbb F_p[x]$. By elementary calculus there is a
polynomial $\lambda_r\in\Bbb Z_{(p)}[x_0,y_0,...,x_r,y_r]$
such that for any $a,b\in R$,
$$D^{(p^r)}(ab)=\lambda_r(D^{(p^0)}\otimes 1,
1\otimes D^{(p^0)},...,D^{(p^r)}\otimes 1,1\otimes D^{(p^r)})(a\otimes b)$$
Using this we explain the left invariant differential operators corresponding to
the Dieudonn\'e elements of the Cartier dual. For example, when $m=0$ we have
$\lambda_r\equiv\phi_r$ (mod $p$), where $\phi_r$ is the polynomial in the definition
of the addition of $\Cal W$. This gives another understanding of $\phi_r$.
This also gives a hint to the understanding via crystalline cohomology, because
differential operator is a background of divided power.  

\vskip 20pt
\flushpar{\bf Acknowledgements.} I wish to thank F. Oort, J.-M. Fontaine and P. Berthelot
for their motivative discussions with me on this topic.

\vskip 30pt
\flushpar {\bf 1. Some prerequisites on the calculus about group schemes}

\vskip 10pt
Thoughout this paper we denote by $k$ a perfect field of characteristic $p>0$.
A finite group scheme $G=\text{Spec}(R)$ over $k$ is called {\sl infinitesimal} if $R$
is a local ring. If $G$ is commutative, denote by $G^D$ its Cartier dual. Denote by
$\goth A\goth b_{\text{inf}}^{\text{inf}}$ the category of finite commutative infinitesimal
group schemes whose Cartier duals are also infinitesimal. An object
$G\in\text{Ob}(\goth A\goth b_{\text{inf}}^{\text{inf}})$ is viewed as an ``additive group
scheme'', and we denote by $a_G:G\times_kG\to G$ the binary operation (addition),
$0_G:\text{Spec}(k)\to G$ the zero section, and $-1_G:G\to G$ the inverse. For convenience
we denote by $a_{G,n}:G\times_k\buildrel{n}\over\cdots\times_kG\to G$ the addition of $n$ factors
($n\ge 2$).

The facts in this section are basically well-known.
We collect them together as an easy reference, but omit the proofs (which can all be found
in [Li2]).

Let $G=\text{Spec}(R)\in\text{Ob}(\goth A\goth b_{\text{inf}}^{\text{inf}})$. A
$k$-differential operator $D:R\to R$ is called {\sl left invariant} if the following
diagram is commutative
$$\CD R @>D>> R \cr @VV{a_G^*}V @VV{a_G^*}V \cr
R\otimes_kR @>{\text{id}_R\otimes_kD}>> R\otimes_kR \endCD\tag 1.1$$
Denote by $\bar D=0_G^*\circ D:R\to k$, then by (1.1) we have
$$D=(\text{id}_R\otimes_k\bar D)\circ a_G^*:R\to R \tag 1.2$$
Hence $D$ is determined by $\bar D$. Denote by $Diff(G/k)$ the set of all
left invariant $k$-differential operators of $G$. Then (1.2) gives a canonical isomorphism
of $k$-linear spaces
$$Diff(G/k)\to Hom_k(R,k)\cong R^D \tag 1.3$$
where $R^D$ is the structure ring of $G^D$. On the other hand, $Diff(G/k)$ has a
$k$-algebra structure, and $Lie(G/k)$ can be viewed as a Lie subalgebra of $Diff(G/k)$.
Note that $Lie(G/k)$ is a ``$p$-Lie algebra'', i.e. for any $\theta\in Lie(G/k)$
we have $\theta^p\in Lie(G/k)$. 

For any $D,D'\in Diff(G/k)$, we have
$$\aligned D\circ D' & = (\text{id}_R\otimes_k\bar D)\circ a_G^*\circ (\text{id}_R\otimes_k\bar D')\circ a_G^* \cr
& = (\text{id}_R\otimes_k\bar D)\circ (\text{id}_{R\otimes_kR}\otimes_k\bar D')\circ (a_G^*\otimes_k
\text{id}_R)\circ a_G^* \cr
& = (\text{id}_R\otimes_k\bar D\otimes_k\bar D')\circ a_{G,3}^*  \endaligned\tag 1.4$$
Hence
$$\aligned \overline{D\circ D'} & =0_G^*\circ (\text{id}_R\otimes_k\bar D\otimes_k\bar D')
\circ a_{G,3}^* \cr
& = (\bar D\otimes_k\bar D')\circ a_G^* \endaligned\tag 1.5$$
Note that the multiplication of $R^D$ is given
by $(\alpha ,\beta )\mapsto (\alpha\otimes_k\beta )\circ a_G^*$ ($\forall\alpha ,\beta\in R^D$).
Thus (1.4) means that the map $Diff(G/k)\to R^D$ in (1.3) is a $k$-algebra isomorphism.

Let $I_0\subset R$ be the maximal ideal, and denote $\omega_{G/k}=I_0/I_0^2$. Then canonically
$$\Omega_{R/k}^1\cong\omega_{G/k}\otimes_kR \tag 1.6$$
and
$$Lie(G/k)\cong\omega_{G/k}^{\vee}=Hom_k(\omega_{G/k},k) \tag 1.7$$
as $k$-linear speces.
Take $x_1,...,x_n\in I_0$ such that their images $\bar x_1,...,\bar x_n\in\omega_{G/k}$ form
a $k$-basis of $\omega_{G/k}$. Then by (1.6) there exist $D_1,...,D_n\in Der_k(R,R)$ such that
$D_i(x_j)=\delta_{ij}$ ($1\le i,j\le n$). From this one sees that all monomials
$x_1^{i_1}\cdots x_n^{i_n}$ ($0\le i_1,...,i_n<p$) are linearly independent over $k$.
In particular, if $F_{G/k}=0$, then $R\cong k[x_1,...,x_n]/(x_1^p,...,x_n^p)$.
By induction (on $m$ such that $F_{G/k}^m=0$) one sees that in general case there exist
a $k$-algebra isomorphism
$$R\cong k[x_1,...,x_n]/(x_1^{p^{i_1}},...,x_n^{p^{i_n}}) \tpt (0<i_1\le ...\le i_n) \tag 1.8$$
Summarizing,

\vskip 10pt
\flushpar{\bf Fact 1.1}. Let $G=\text{Spec}(R)\in\text{Ob}(\goth A\goth b_{\text{inf}}^{\text{inf}})$.
Let $I_0\subset R$ be the maximal ideal, and denote $\omega_{G/k}=I_0/I_0^2$. 

\vskip 5pt\leftskip=30pt\parindent=-10pt
i) The $k$-algebra $Diff(G/k)$ of left invariant differential operators is canonically isomorphic
to $R^D$ (the structure ring of $G^D$).

ii) The $k$-Lie algebra $Lie(G/k)$ of left invariant derivations is a $p$-Lie subalgebra
of $Diff(G/k)$ over $k$, and there are canonical isomorphisms (1.6) and (1.7).

iii) Let $n=\dim_k(\omega_{G/k})=\dim_k(Lie(G/k))$. Then $R$ has structure (1.8)
as a $k$-algebra, where $x_1,...,x_n$ generate the maximal ideal of $R$.

\vskip 10pt\leftskip=0pt\parindent=20pt
Denote by $\Bbb G_{a/k}=\text{Spec}[t]$ whose (additive) group scheme structure is given
by $a^*(t)=t\otimes_k1+1\otimes_kt$, $0^*(t)=0$, $(-1)^*(t)=-t$. It is easy to see that
any $k$-group scheme homomorphism $f:G\to\Bbb G_{a/k}$ is uniquely determined by an
element $x=f^*(t)\in I_0$ satisfying $a_G^*(x)=x\otimes_k1+1\otimes_kx$. Let
$$\alpha (G)=\{ x\in R|a^*(x)=x\otimes_k1+1\otimes_kx\} \tag 1.9$$
Obviously $\alpha (G)$ has a $k$-linear space structure, called the {\sl $\alpha$-module} of $G$.
The above argument shows that there is a canonical one-to-one correspondence
$$Hom_k(G,\Bbb G_{a/k})\leftrightarrow\alpha (G) \tag 1.10$$
View $x$ as a linear functional on $R^D$, then (1.9) means for any $y,y'\in R^D$ we have
$x(yy')=x(y)y'+yx(y')$, i.e. $x\in Der_k(R^D,R^D)$. Thus by Fact 1.1 there is a canonical
one-to-one correspondence
$$\alpha (G)\leftrightarrow Lie(G^D/k) \tag 1.11$$
Furthermore, for any homomorphism $f:G\to G'$ in $\goth A\goth b_{\text{inf}}^{\text{inf}}$,
$H=\text{coker}(f)$ gives an exact sequence
$0\to H^D\to G^{\prime D}\buildrel{f^D}\over\to G^D$, which implies a left exact sequence
$$0\to Lie (H^D/k)\to Lie(G^{\prime D}/k)\buildrel{f^D_*}\over\longrightarrow Lie(G^D/k) \tag 1.12$$
By (1.11), this gives a canonical left exact sequence
$$0\to\alpha (H)\to\alpha (G')\buildrel{f^*}\over\longrightarrow\alpha (G) \tag 1.13$$
Since $G^D\in\text{Ob}(\goth A\goth b_{\text{inf}}^{\text{inf}})$, when $G\ne 0$ we have
$\omega_{G^D/k}\ne 0$, hence by (1.7) we have $Lie(G^D/k)\ne 0$, therefore by (1.11) we
have $\alpha (G)\ne 0$. Let $f:G\to\Bbb G_{a/k}$ be a non-zero homomorphism, then
$\text{im}(f)\cong\alpha_{p^n}$ for some $n>0$, hence $f$ induces an epimorphism
$f':G\to\alpha_p$. By Fact 1.1.iii) we see $\alpha_p$ is a simple object in
$\goth A\goth b_{\text{inf}}^{\text{inf}}$. Denote by $l(G)$ the length of $G$
in $\goth A\goth b_{\text{inf}}^{\text{inf}}$.
Note that $\text{ker}(f')\in\text{Ob}(\goth A\goth b_{\text{inf}}^{\text{inf}})$.
By induction on $l(G)$ we get a filtration
$$G=G_0\supset G_1\supset ...\supset G_n=0 \tag 1.14$$
whose factors are all isomorphic to $\alpha_p$. Hence
$$n=l(G),\ \ \dim_k(R)=p^{l(G)} \tag 1.15$$

\vskip 10pt
\flushpar{\bf Example 1.1}. Let $G=\alpha_{p^n}\cong\Bbb G_{a/k}[F_{G/k}^n]$. Then
$R\cong k[x]/(x^{p^n})$ (where $x$ is the image of $t$, as above), and $\alpha (G)$ has a $k$-basis
$x,x^p,...,x^{p^{n-1}}$.

\vskip 10pt
For any $x\in\alpha (G)$, obviously $x^p\in\alpha (G)$. If there are $x_1,...,x_n\in\alpha (G)$
such that (1.8) holds, then by Example 1.1 one sees that
$$G\cong\alpha_{p^{i_1}}\times_k\cdots\times_k\alpha_{p^{i_n}} \tag 1.16$$

If $F_{G^D/k}=0$, then by Fact 1.1.iii), (1.15) and (1.10) we see that
$$l(G)=l(G^D)=\dim_k(\omega_{G^D/k})=\dim_k(\alpha (G)) \tag 1.17$$
Obviously $\alpha (\alpha_p)\cong k$. Hence by (1.13) and induction we see that $R$
is generated by $\alpha (G)$ as a $k$-algebra, hence (1.16) holds. If $F_{G/k}=0$ also,
then $i_1=...=i_n=1$, i.e. $G\cong\alpha_p^n$. In this case we say $G$ is an {\sl $\alpha$-group}.
The above argument shows that $G$ is an $\alpha$-group if and only if both
$F_{G/k}=0$ and $V_{G/k}=0$.

It is easy to see that a homomorphism $\alpha_p^m\to\alpha_p^n$ is equivalent to a $k$-linear
homomorphism $k^n\cong\alpha (\alpha_p^n)\to\alpha (\alpha_p^m)\cong k^m$. For any
$G\in\text{Ob}(\goth A\goth b_{\text{inf}}^{\text{inf}})$, there is a largest
$\alpha$-subgroup $H\subset G$, in fact
$$H=\text{ker}(F_{G/k})\cap\text{ker}(V_{G^{(p^{-1})}/k}) \tag 1.18$$
Note that any homomorphism $\alpha_p\to G$ factors through $H$, hence the above argument
shows that $Hom_k(\alpha_p,G)$ has naturally a $k$-linear space structure, and
$$\dim_k(Hom_k(\alpha_p,G))=\dim_k(\alpha (H)) \tag 1.19$$
We call $\dim_k(Hom_k(\alpha_p,G))$ the {\sl $a$-number} of $G$, denoted by $a(G)$.
The filtration (1.14) shows that $a(G)>0$ if $G\ne 0$. Summarizing,

\vskip 10pt
\flushpar{\bf Fact 1.2}. Let $G=\text{Spec}(R)\in\text{Ob}(\goth A\goth b_{\text{inf}}^{\text{inf}})$.

\vskip 5pt\leftskip=30pt\parindent=-10pt
i) There are canonical one-to-one correspondences
$$Hom_k(G,\Bbb G_{a/k})\leftrightarrow\alpha (G)\cong Lie(G^D/k) \tag 1.20$$
Furthermore, for any homomorphism $f:G\to G'$ in $\goth A\goth b_{\text{inf}}^{\text{inf}}$,
there is a canonical left exact sequence
$$0\to\alpha (\text{coker}(f))\to\alpha (G')\buildrel{f^*}\over\longrightarrow\alpha (G) \tag 1.21$$

ii) $G$ has a largest $\alpha$-subgroup $H$ given by (1.18), and (1.19) holds.
In particular $a(G)\ne 0$ when $G\ne 0$.

iii) There is a filtration (1.14) whose factors are all isomorphic to $\alpha_p$, where
$n=l(G)$. Furthermore $\dim_k(R)=p^{l(G)}$.

iv) If $V_{G/k}=0$, then $R$
is generated by $\alpha (G)$ as a $k$-algebra, and $G$ has a decomposition (1.16).

v) $G$ is an $\alpha$-group if and only if both $F_{G/k}=0$ and $V_{G/k}=0$.

\vskip 10pt\leftskip=0pt\parindent=20pt
Let $H\cong\alpha_p^r\subset G$ be as in Fact 1.2.ii), where $r=a(G)$.
For any $H'\cong\alpha_p\subset H$, one can take $H^{\prime\prime}\subset H$ such that
$H=H'\times_kH^{\prime\prime}$. If $H^{\prime\prime}\ne 0$,
let $G_1=H/H^{\prime\prime}$ and $H_1$ be the
largest $\alpha$-subgroup of $G_1$, then $H'$ can be viewed as a subgroup scheme of
$H_1$. If $H_1\ne H'$, replace $G$ and $H$ by $G_1$ and $H_1$ respectively and repeat
the above argument. By induction one eventually gets a quotient group scheme $G'$ of
$G$ such that $H'\to G'$ is injective and $a(G')=1$. Since $H'$ is arbitrary, we can get
$G_1,...,G_r\in\text{Ob}(\goth A\goth b_{\text{inf}}^{\text{inf}})$ with
$a(G_1)=...=a(G_r)=1$ such that there is a homomorphism
$$f:G\to G_1\times_k\cdots\times_kG_r \tag 1.22$$
whose restriction on $H$ is injective. This means $f$ is injective by Fact 1.2.ii), because
$\text{ker}(f)\cap H=0$ (hence $a(\text{ker}(f))=0$). Thus we have

\vskip 10pt
\flushpar{\bf Corollary 1.3}. For any
$G\in\text{Ob}(\goth A\goth b_{\text{inf}}^{\text{inf}})$ with $a(G)=r$,
there exist $G_1,...,G_r\in\text{Ob}(\goth A\goth b_{\text{inf}}^{\text{inf}})$ with
$a(G_1)=...=a(G_r)=1$ such that there is a monomorphism $f$ as in (1.22).

\vskip 30pt
\flushpar {\bf 2. Some prerequisites about the Witt scheme}

\vskip 10pt
The definition of Witt scheme is as follows (see [Mu1, \S 23]). Fix a prime number $p$.
For each integer $n\ge 0$, let
$$w_n(x_0,..., x_n)=\sum_{i=0}^{n}p^ix_i^{p^{n-i}}\in\Bbb Z[x_0,..., x_n]
\tpt (n=0,1,...) \tag 2.1$$
One can inductively define $\phi_n\in\Bbb Q[x_0,y_0,...,x_n,y_n]$ ($n=0,1,2,...$)
such that for any $n\ge 0$,
$$w_n(x_0,..., x_n)+w_n(y_0,..., y_n)=
w_n(\phi_0(x_0, y_0),...,\phi_n(x_0, y_0,..., x_n, y_n)) \tag 2.2$$
It is easy to see two facts: One is that for any $f(t_1,...,t_m)\in\Bbb Z[t_1,...,t_m]$,
the difference $f(t_1,...,t_m)^p-f(t_1^p,...,t_m^p)$ is divisible by $p$; the other is that
$(x+py)^{p^m}-x^{p^m}$ is divisible by $p^{m+1}$. Hence the difference of
$$w_n(\phi_0(x_0, y_0)^p,...,\phi_n(x_0, y_0,..., x_n, y_n)^p) \tag 2.3$$
and
$$w_n(\phi_0(x_0^p, y_0^p),...,\phi_n(x_0^p, y_0^p,..., x_n^p, y_n^p))
\tag 2.4$$
is divisible by $p^{n+1}$. But this difference is equal to $p^{n+1}(x_{n+1}+y_{n+1})$,
hence $\phi_n\in\Bbb Z[x_0,y_0,...,x_n,y_n]$. Similarly there are
$\psi_n\in\Bbb Z[x_0,y_0,...,x_n,y_n]$ ($n=0,1,2,...$) such that for any $n\ge 0$,
$$w_n(x_0,..., x_n)w_n(y_0,..., y_n)=
w_n(\psi_0(x_0, y_0),...,\psi_n(x_0, y_0,..., x_n, y_n)) \tag 2.5$$

\vskip 10pt
\flushpar{\bf Fact 2.1}. Let $A_n=\Bbb Z[x_0,..., x_n]$,
$\Cal W_n=\text{Spec}A_n$. Define ring homomorphisms
$a^*:A_n\to A_n\otimes A_n\cong\Bbb Z[x_0, ...,x_n, y_0,..., y_n]$
and $m^*:A_n\to A_n\otimes A_n$ by  
$a^*(x_i)=\phi_i(x_0,y_0,...,x_i,y_i)$ and $m^*(x_i)=\psi_i(x_0,y_0,...,x_i,y_i)$
respectively. Then $\Cal W_n$ has a commutative ring scheme
structure (with 1), whose addition $a$ (resp. multiplication $m$) is given by $a^*$
(resp. $m^*$).
 
\vskip 10pt
\flushpar {\sl Proof}. Let $B=A_n\otimes\Bbb Q\cong\Bbb Q[x_0, ...,x_n]$.
Note that $w_0, ...,w_n$ give a ring homomorphism $f:A_n\to A_n$ such that
$f_{\Bbb Q}=f\otimes\text{id}_{\Bbb Q}:B\to B$ is an isomorphism.
Define ring homomorphisms $a_B^*:B\to B\otimes B\cong\Bbb Q[x_0, x_1,..., y_0, y_1,...]$
and $m_B^*:B\to B\otimes B$ by $a_B^*(x_n)=x_n+y_n$ and $m_B^*(x_n)=x_ny_n$
respectively. It is obvious that $a_B^*$ and $m_B^*$ define a ring scheme structure
on $\text{Spec}B$.
Denote by $a_{\Bbb Q}^*=a^*\otimes\text{id}_{\Bbb Q}:B\to B\otimes B$, and
$m_{\Bbb Q}^*=m^*\otimes\text{id}_{\Bbb Q}:B\to B\otimes B$. Then by the definitions above
we have $a_B^*\circ f_{\Bbb Q}=(f_{\Bbb Q}\otimes f_{\Bbb Q})\circ a_{\Bbb Q}^*$ and
$m_B^*\circ f_{\Bbb Q}=(f_{\Bbb Q}\otimes f_{\Bbb Q})\circ m_{\Bbb Q}^*$. This shows that $a_{\Bbb Q}^*$ and
$m_{\Bbb Q}^*$ define a commutative ring scheme structure on $\text{Spec}B$.
Since $A_n$ is a subring of $B$, this shows that the addition $a$ and multiplication $m$
satisfy the axioms of ring scheme. \ \ \ Q.E.D.

\vskip 10pt
We call $\Cal W_n$ the {\sl truncated Witt scheme} of level $n$.
Note that $\Cal W_0$ is the ring scheme structure of $\Bbb A^1$.
For any $n>0$ denote by $q_n:\Cal W_n\to\Cal W_{n-1}$ the projection $(x_0,...,x_n)\to (x_0,...,x_{n-1})$,
which is obviously a faithfully flat ring scheme homomorphism. Denote by $\Cal W$ the formal
limit $\varprojlim\limits_n\Cal W_n$, called the {\sl Witt scheme}.
By definition one also sees that the ideal $(x_0)\subset A_{n+1}$ defines a closed subgroup
scheme of $\Cal W_{n+1}$ which is isomorphic to $\Cal W_n$, i.e. the closed immersion
$i_n:\Cal W_n\to\Cal W_{n+1}$ defined by $i_n^*(x_0)=0$, $i_n^*(x_{i+1})=x_i$ ($0\le i\le n$)
is a group scheme homomorphism. Furthermore, the morphism $V=V_n:\Cal W_n\to\Cal W_n$ defined by
$V^*(x_i)=x_{i-1}$ ($\forall i>0$), $V^*(x_0)=0$ (i.e. $(x_0,...,x_n)\to (0,x_0,...,x_{n-1})$)
is an endomomorphism of the additive group scheme structure of $\Cal W_n$, called the {\sl Verschiebung}
of $\Cal W_n$. It is easy to see that
$$\text{ker}(V_n)\cong\text{coker}(V_n)\cong\Cal W_0\cong\Bbb G_{a/\Bbb Z} \tag 2.6$$
The morphism $\tau =\tau_n:\Cal W_0\to\Cal W_n$ defined by
$x_0\mapsto (x_0,0,...,0)$ commutes with the multiplications (but not the additions when $n>0$),
called the {\sl Teichm\"uller lifting}, which is a section of the projection $\Cal W_n\to\Cal W_0$.

From these one sees the following properties.

\vskip 10pt
\flushpar{\bf Corollary 2.1}. Notation as above. For simplicity denote
$X_i=(x_0,...,x_i)$, $Y_i=(y_0,...,y_i)$, $Z_i=(z_0,...,z_i)$ and denote
$\phi_i(x_0,y_0,...,x_i,y_i)$ by $\phi_i(X_i,Y_i)$.

\vskip 5pt\leftskip=30pt\parindent=-10pt
i) For each $i$,
$$\phi_i(X_i,Y_i)-x_i-y_i\in\Bbb Z[x_0,y_0,...,x_{i-1},y_{i-1}] \tag 2.7$$
in particular $\phi_0(x_0,y_0)=x_0+y_0$.

ii) For each $i>0$,
$$\phi_i(0,0,x_0,y_0,...,x_{i-1},y_{i-1})=\phi_{i-1}(x_0,y_0,...,x_{i-1},y_{i-1}) \tag 2.8$$

iii) For each $i$ there hold the co-commutativity
$$\phi_i(X_i,Y_i)=\phi_i(Y_i,X_i) \tag 2.9$$
and co-associativity
$$\aligned & \phi_i(x_0,\phi_0(Y_0,Z_0),x_1,\phi_1(Y_1,Z_1),...,x_i,\phi_i(Y_i,Z_i)) \cr
= & \phi_i(\phi_0(X_0,Y_0),z_1,\phi_1(X_1,Y_1),z_2,...,\phi_i(X_i,Y_i),z_i) \endaligned\tag 2.10$$

iv) For each monomial $cx_0^{m_0}y_0^{m_0^{\prime}}\cdots x_i^{m_i}y_i^{m_i^{\prime}}$ ($c\ne 0$) of
$\phi_i$, it holds that
$$\sum\limits_{j=0}^i(m_j+m_j^{\prime})p^j=p^i \tag 2.11$$
In particular, when $i>0$ we have $m_0+m_0^{\prime}\equiv 0$(mod $p$); when $i>1$ we have
$m_0+m_0^{\prime}+p(m_1+m_1^{\prime})\equiv 0$ (mod $p^2$), etc..

\vskip 10pt\leftskip=0pt\parindent=20pt
Let $k\supset\Bbb F_p$ be a perfect field. Then by Fact 2.1 we see $W(k)=\Cal W(k)$ is a complete
DVR with maximal ideal $(p)$ such that $W(k)/(p)\cong k$. It is elementary that such a DVR is
unique up to isomorphism (whenever it exists). Furthermore, the Teichm\"uller liftings $\tau_n$
($\forall n$)
give a multiplicative map $\tau =\tau (k):k\to W(k)$ which is a section of the projection
$W(k)\to k$ (also called ``Teichm\"uller lifting''). Any element in $W(k)$ can be uniquely
written as a series
$$\sum_{i=0}^{\infty}p^i\tau (a_i^{p^{-i}})
=\tau (a_0)+p\tau (a_1^{p^{-1}})+p^2\tau (a_2^{p^{-2}})+\cdots \tpt (a_i\in k) \tag 2.12$$
which can be simply expressed as $(a_0,a_1,a_2,...)$, called a {\sl Witt vector}. Note that
the unit element in $W(k)$ is $(1,0,0,...)$. By above, the
addition and multiplication of $W(k)$ can be given by
$$\aligned & (a_0,a_1,a_2,...)+(b_0,b_1,b_2,...) \cr
= & (\phi_0(a_0,b_0),\phi_1(a_0,b_0,a_1,b_1),
\phi_2(a_0,b_0,a_1,b_1,a_2,b_2),...) \endaligned\tag 2.13$$
and
$$\aligned & (a_0,a_1,a_2,...)\cdot (b_0,b_1,b_2,...) \cr
= & (\psi_0(a_0,b_0),\psi_1(a_0,b_0,a_1,b_1),
\psi_2(a_0,b_0,a_1,b_1,a_2,b_2),...) \endaligned\tag 2.14$$
respectively.

By definition we have $p\cdot (a_0,a_1,a_2,...)=(0,a_0^p,a_1^p,a_2^p,...)$, which is equal to
$F_{\Cal W\otimes k/k}\circ V((a_0,a_1,a_2,...))=V\circ F_{\Cal W\otimes k/k}((a_0,a_1,a_2,...))$.
Since each $\Cal W_n$ is a variety, this shows that
$$F_{\Cal W\otimes k/k}\circ V=V\circ F_{\Cal W\otimes k/k}=p_{\Cal W\otimes k}:
\Cal W\otimes k\to\Cal W\otimes k \tag 2.15$$
Note that $(a_0,a_1,a_2,...)\mapsto F_{\Cal W\otimes k/k}((a_0,a_1,a_2,...))=(a_0^p,a_1^p,a_2^p,...)$
is a ring automorphism of $W(k)$, called the ``Frobenius'' of $W(k)$ and denoted by $\sigma$.
Summarizing we have

\vskip 10pt
\flushpar{\bf Corollary 2.1}. Let $k\supset\Bbb F_p$ be a perfect field.

\vskip 5pt\leftskip=30pt\parindent=-10pt
i) Up to isomorphism, $W(k)=\Cal W(k)$ is the unique complete
DVR with maximal ideal $(p)$ such that $W(k)/(p)\cong k$. The Teichm\"uller liftings $\tau_n$
($\forall n$) give a multiplicative map (Teichm\"uller lifting)
$\tau =\tau (k):k\to W(k)$ which is a section of the projection $W(k)\to k$.

ii) Any element in $W(k)$ can be uniquely written as a series of the form (2.12), which
can be simply expressed as a Witt vector $(a_0,a_1,a_2,...)$. In this way the
addition and multiplication of $W(k)$ can be given by (2.13) and (2.14)
respectively, with unit element $(1,0,0,...)$.

iii) (2.15) holds on $\Cal W\otimes k$.

iv) The Frobenius $\sigma :(a_0,a_1,a_2,...)\mapsto (a_0^p,a_1^p,a_2^p,...)$ is a ring
automorphism of $W(k)$, and $\sigma\circ V=V\circ\sigma =p\cdot :W(k)\to W(k)$, where
$V$ satisfies $V((a_0,a_1,a_2,...))=(0,a_0,a_1,a_2,...)$.

\vskip 10pt\leftskip=0pt\parindent=20pt
Now let $k=\Bbb F_p$. Denote by $W_n=\Cal W_n\otimes k$, and $W_{n,m}=W_n[F_{W_n/k}^{m+1}]$
for any $m\ge 0$. Then each $W_{n,m}$ is a finite commutative group scheme,
$W_{0,0}\cong\alpha_p$, and $W_{n,m}$ can be viewed as a closed subgroup scheme of
$W_{n+1,m}$ or $W_{n,m+1}$. As a scheme,
$$W_{n,m}\cong\Bbb F_p[x_0,...,x_n]/(x_0^{p^{m+1}},...,x_n^{p^{m+1}}) \tag 2.16$$

\vskip 10pt
\flushpar{\bf Fact 2.2}. For any $m\ge 0$, the group scheme endomorphism $V_{n,m}$ of $W_{n,m}$
induced by $V_n$ is equal to $V_{W_{n,m}/\Bbb F_p}$.

\vskip 10pt
\flushpar {\sl Proof}. Note that $F_{W_{n,m+1}/\Bbb F_p}$ factors through $W_{n,m}$, inducing a
faithfully flat homomorphism $q_m:W_{n,m+1}\to W_{n,m}$, whose composition with the embedding
$i_m:W_{n,m}\to W_{n,m+1}$ is equal to $F_{W_{n,m+1}/\Bbb F_p}$. By (2.15) we have
$$V_{n,m+1}\circ F_{W_{n,m+1}/\Bbb F_p}=p_{W_{n,m+1}} \tag 2.17$$
By the commutative diagram
$$\CD W_{n,m+1} @>{q_m}>> W_{n,m} @>{V_{n,m}}>> W_{n,m} \cr
@VV{\text{id}}V @VV{i_m}V @VV{i_m}V \cr
W_{n,m+1} @>{F_{W_{n,m+1}/\Bbb F_p}}>> W_{n,m+1} @>{V_{n,m+1}}>> W_{n,m+1} \endCD\tag 2.18$$
we see that
$$i_m\circ V_{n,m}\circ q_m=p_{W_{n,m+1}} \tag 2.19$$
On the other hand, by the theory of commutative group schemes (cf. [Dem]) we have
$$V_{W_{n,m}/\Bbb F_p}\circ F_{W_{n,m}/\Bbb F_p}=p_{W_{n,m}} \tag 2.20$$
Hence similarly we have
$$i_m\circ V_{W_{n,m}/\Bbb F_p}\circ q_m=p_{W_{n,m+1}} \tag 2.21$$
Comparing (2.19) and (2.21), noting that $i_m$ is injective and $q_m$ is surjective, we
see that $V_{n,m}=V_{W_{n,m}/\Bbb F_p}$. \ \ \ Q.E.D.

\vskip 10pt
By (2.6) and Fact 2.2 we see that
$$W_n[V]\cong\text{coker}(V_{W_n/k})\cong\Bbb G_{a/k} \tag 2.22$$
Hence $a(W_n)=1$, and for any $m\ge 0$,
$$W_{n,m}[V]\cong\text{coker}(V_{W_{n,m}/k})\cong\alpha_{p^{m+1}}  \tag 2.23$$
Thus $a(W_{n,m})=a(W_{n,m}^D)=1$. Summarizing we have

\vskip 10pt
\flushpar{\bf Corollary 2.3}. Let $k=\Bbb F_p$, and denote by $W_n=\Cal W_n\otimes k$ and
$W_{n,m}=W_n[F_{W_n/k}^{m+1}]$ for any $m\ge 0$. Then

\vskip 5pt\leftskip=30pt\parindent=-10pt
i) Each $W_{n,m}\in\text{Ob}(\goth A\goth b_{\text{inf}}^{\text{inf}})$,
$W_{0,0}\cong\alpha_p$, and $W_{n,m}$ can be viewed as a closed subgroup scheme of
$W_{n+1,m}$ or $W_{n,m+1}$. The scheme structure of $W_{n,m}$ is given by (2.16).

ii) (2.22) and (2.23) hold.

iii) For any $m,n$, $a(W_n)=a(W_{n,m})=a(W_{n,m}^D)=1$.

\vskip 30pt\leftskip=0pt\parindent=20pt
\flushpar {\bf 3. Dieudonn\'e elements and Dieudonn\'e modules}

\vskip 10pt
Let $k\supset\Bbb F_p$ be a perfect field. We use the notation as above.
Note that for any finite group scheme $G=\text{Spec}R$ over $k$ and any $m\in\Bbb Z$,
$G^{(p^m)}=\text{Spec}R^{(p^m)}$ makes sense since $k$ is perfect. If $G$
is commutative, any element $x\in R$ can be viewed as an element of $R^{(p^m)}$
via the ring isomorphism $R^{(p^m)}\cong R$, which will also be denoted by $x$,
by abuse of notation. Hence $V_{G/k}^{m*}x$ can also be viewed as an element of $R$,
simply denoted as $V^{m*}x$.

\vskip 10pt
\flushpar{\bf Definition 3.1}. Let $G=\text{Spec}R$ be a finite commutative (additive) group scheme
such that $V^{n+1}_{G/k}=0$. Denote by $M=\text{ker}(0_G^*)\subset R$. An element
$x\in M$ is called a {\sl Dieudonn\'e element} if
$$a_G^*(x)=\phi_n(V^{n*}x\otimes_k1,1\otimes_kV^{n*}x,V^{n-1*}x\otimes_k1,
1\otimes_kV^{n-1*}x,...,x\otimes_k1,1\otimes_kx) \tag 3.1$$
and called a {\sl Dieudonn\'e generator} if in addition that $R$ is generated by
$x$, $V^*x$, $V^{2*}x$,... $V^{n*}x$ as a $k$-algebra.

\vskip 10pt
Note that by Corollary 2.1.ii), the above definition is independent of the choice of $n$.

In the following denote by $W_n=\Cal W_n\otimes k$, and $W_{n,m}=\Cal W_{n,m}\otimes k$.
If $f:G\to W_n$ is a homomorphism of group schemes over $k$, then by Fact 2.1 we see
that $f^*(x_n)$ is a Dieudonn\'e element of $G$; conversely, if $x$ is a
Dieudonn\'e element of $G=\text{Spec}R$, then one can define a $k$-morphism
$f:G\to W_n$ by $f^*(x_i)=V^{n-i*}_{G/k}(x)$ ($0\le i\le n$), which is a group
scheme homomorphism by Fact 2.1 and Definition 3.1. Thus we can denote $f=f_x$.

Let  $G=\text{Spec}R$ be a finite commutative group scheme
such that $V^{n+1}_{G/k}=0$.
Denote by $D(G)\subset R$ the set of all Dieudonn\'e elements. For any $x,y\in D(G)$,
$$\CD f_x+f_y:G @>{\Delta}>> G\times_kG @>{(f_x,f_y)}>> W_n\times W_n
@>{a}>> W_n \endCD\tag 3.2$$
is also a homomorphism of $k$-group schemes, hence $(f_x+f_y)^*(x_n)\in D(G)$, denoted by
$x\dotplus y$. By Fact 2.1 we have
$$x\dotplus y=\phi_n(V^{n*}x,V^{n*}y,V^{n-1*}x,V^{n-1*}y,...,x,y) \tag 3.3$$
Furthermore we have $-x=-1_G^*(x)\in D(G)$ and $0\in D(G)$. By Corollary 2.1.iii)
we have
$$\aligned & x\dotplus y=y\dotplus x,\ (x\dotplus y)\dotplus z=x\dotplus (y\dotplus z), \cr
& x\dotplus 0=x,\ x\dotplus (-x)=0\ \ (\forall x,y,z\in D(G)) \endaligned\tag 3.4$$
This means $D(G)$ has an ``additive'' commutative group structure under $\dotplus$.

Note that an element $c\in W_n(k)$ can be viewed as a $k$-morphism $\text{Spec}(k)\to W_n$.
Any $x\in D(G)$ and $c\in W_n(k)$ gives a $k$-group scheme homomorphism
$$\CD G\cong\text{Spec}(k)\times_kG @>{(c,f_x)}>> W_n\times W_n
@>{m}>> W_n \endCD\tag 3.5$$
This corresponds to a Dieudonn\'e element, denoted by $c\dot\times x$. If
$c=(c_0,c_1,...,c_n)$, by the definition of Witt scheme it is easy to see that
$$c\dot\times x=c_0^{p^n}x\dotplus c_1^{p^{n-1}}p^*x\dotplus\cdots\dotplus c_np^{n*}x \tag 3.6$$
In particular, if $c_1=...=c_n=0$, then $c\dot\times x=c_0^{p^n}x$. Thus for any $c\in k$
and $x\in D(G)$ we have $cx\in D(G)$.

It is easy to see that if $x\in D(G)$ then $V^*x\in D(G)$ and $x^p=F^*x\in D(G)$,
and by (3.5) it is easy to see that for any $c\in W_n(k)$,
$$F^*(c\dot\times x)=c^{\sigma}\dot\times F^*x,\ \  V^*(c\dot\times x)=c^{\sigma^{-1}}\dot\times V^*x
\tag 3.7$$
Let $W=W(k)=\varprojlim\limits_nW_n(k)$.
Define the following $W$-algebra (usually non-commutative)
$$A=W[F,V]/(FV-p,VF-p,Fa-a^\sigma F,Va-a^{\sigma^{-1}}V\ \forall a\in W) \tag 3.8$$
Then $D(G)$ has an $A$-mudule structure, thus we call $D(G)$ the {\sl Dieudonn\'e module} of $G$.

Obviously $\alpha (G)\subset D(G)$, and $D(G)=\alpha (G)$ when $V_{G/k}=0$. In particular,
when $G$ is an $\alpha$-group we have $\alpha (G)\cong k^{\oplus r}$ ($r=a(G)$), where each
copy of $k$ is viewed as the $A$-module $A/A(F,V)$.

It is easy to see that for a homomorphism $f:G=\text{Spec}R\to G'=\text{Spec}R'$ in
$\goth A\goth b_{\text{inf}}^{\text{inf}}$, the $k$-algebra homomorphism $f^*:R'\to R$
satisfies $f^*(D(G'))\subset D(G)$, hence induces an $A$-module homomorphism
$D(G')\to D(G)$, also denoted by $f^*$. If $f$ is an epimorphism and $H=\text{ker}(f)$, then
for any $x\in\text{ker}(D(G)\to D(H))$, obviously the restriction of $f_x:G\to W_n$ on $H$
is 0, hence $f_x$ factors through $G'$, thus $x\in f^*(D(G'))$; on the other hand,
obviously $f^*(D(G'))\subset\text{ker}(D(G)\to D(H))$,
hence $f^*(D(G'))=\text{ker}(D(G)\to D(H))$.

In particular, if $H\cong\alpha_p$, then there is an exact sequence
$0\to D(G/H)\to D(G)\to D(H)$ of $A$-modules, hence by Fact 1.2.ii) and induction
we see that $D(G)$ has finite length as a $W$-module. Summarizing we have

\vskip 10pt
\flushpar{\bf Proposition 3.1}. Let $k$ be a perfect field of characteristic $p>0$,
and $G=\text{Spec}R$ be an infinitesimal commutative group scheme such that $V^{n+1*}_{G/k}=0$.
Let $D(G)\subset R$ be the set of all Dieudonn\'e elements of $G$. Then

\vskip 5pt\leftskip=30pt\parindent=-10pt
i) Each element in $D(G)$ is equivalent to a $k$-group scheme homomorphism $f_x:G\to W_n$
(satisfying $f_x^*(x_n)=x$).

ii) $D(G)$ has an $A$-mudule structure: the sum of any $x,y\in D(G)$ is given by
$x\dotplus y$ in (3.3), the zero element is 0, and the negative element of any $x\in D(G)$ is $-x$;
for any $c=(c_0,c_1,...,c_n)\in W_n(k)$ and $x\in D(G)$, the product of $c$ and $x$ is given by
$c\dot\times x$ in (3.6); $Fx=F_{G/k}^*x=x^p$, $Vx=V_{G/k}^*x$.
Furthermore, for any $c\in k$ and $x\in D(G)$ we have $cx\in D(G)$, which is equal to
$\tau (c^{p^{-n}})\dot\times x$. As a $W$-module $D(G)$ has finite length.

iii) Let $f:G=\text{Spec}R\to G'=\text{Spec}R'$ be a homomorphism in
$\goth A\goth b_{\text{inf}}^{\text{inf}}$. Then $k$-algebra homomorphism $f^*:R'\to R$
induces an $A$-module homomorphism $f^*:D(G')\to D(G)$ canonically. Hence $D$ can be
viewed as a contravariant functor from $\goth A\goth b_{\text{inf}}^{\text{inf}}$ to
$\goth M_A$ (the category of $A$-modules). Furthermore, if $f$ is an epimorphism,
then the following sequence is exact:
$$\CD 0\to D(G') @>{f^*}>> D(G) \to D(\text{ker}(f)) \endCD\tag 3.9$$

iv) $\alpha (G)\subset D(G)$, and $D(G)=\alpha (G)$ when $V_{G/k}=0$. In particular,
if $G$ is an $\alpha$-group, then
$D(G)\cong (A/A(F,V))^{\oplus r}$ ($r=a(G)$).

v) $\Cal W_n$ has the following universality: for any perfect field $k\supset\Bbb F_p$,
any $G=\text{Spec}(A)\in\text{Ob}(\goth A\goth b_{\text{inf}}^{\text{inf}})$
and any $x\in D(G)$ satisfying $V_{G/k}^{n+1*}x=0$, there exists a unique $k$-group scheme
homomorphism $f:G\to\Cal W_n\otimes k$
such that $f^*(x_n)=x$. And $\Cal W_{n,m}=(\Cal W_n\otimes\Bbb F_p)[F^{m+1}]$ has the following
universality: for any perfect field
$k\supset\Bbb F_p$, any $G=\text{Spec}(A)\in\text{Ob}(\goth A\goth b_{\text{inf}}^{\text{inf}})$
and any $x\in D(G)$ satisfying $V_{G/k}^{n+1*}x=0$ and $F_{G/k}^{m+1*}x=0$, there exists a unique
$k$-group scheme homomorphism $f:G\to\Cal W_{n,m}\otimes k$
such that $f^*(x_n)=x$.

\vskip 10pt\leftskip=0pt\parindent=20pt
Denote by $W_{n,m}=\Cal W_{n,m}\otimes k$.
From Proposition 3.1 one easily sees that if $G=\text{Spec}R$ is a closed subgroup scheme of
some $W_{n,m}^r$, then $R$ is generated by $D(G)$ as a $k$-algebra, and the maximal ideal
$M=\text{ker}(0^*)\subset R$ is generated by $D(G)$; conversely, if $R$ is generated by $D(G)$ as
a $k$-algebra, then $G$ can be embedded into some $W_{n,m}^r$ as a closed subgroup scheme over $k$.
The following is the ``key fact''.

\vskip 10pt
\flushpar{\bf Proposition 3.2}. Let $k$ be a perfect field of characteristic $p>0$,
and $G=\text{Spec}(R)\in\text{Ob}(\goth A\goth b_{\text{inf}}^{\text{inf}})$. Then
$G$ can be embedded into some $W_{n,m}^r$ as a closed subgroup scheme over $k$.

\vskip 10pt
\flushpar {\sl Proof}. By Corollary 1.3 it can be reduced to the case $a(G)=1$.
For simplicity assume $F_{G/k}^{m+1}=0$ and $V_{G/k}^{n+1}=0$. By Fact 1.2.ii) we can
take a closed subgroup scheme $H\subset G$ such that $G/H\cong\alpha_p$. Note that $a(H)\le 1$,
hence by induction we need only to prove that if $H$ is a closed subgroup scheme of $W_{n,m-1}$,
then $G$ can be embedded into $W_{n,m}$ as a closed subgroup scheme.
Note that $G$ has a closed subgroup scheme $H'=G[F]\cap G[V]\cong\alpha_p$, and a homomorphism
$f:G\to W_{n,m}$ is injective if and only if its restriction on
$H'$ is injective, because if $\text{ker}(f)\cap H'=0$ then
$a(\text{ker}(f))=0$, hence $\text{ker}(f)=0$ by Fact 1.2.i).

For times we will use the fact that the Cartier dual is an anti-equivalence of
$\goth A\goth b_{\text{inf}}^{\text{inf}}$.

Let $G_0$ be the push-out of $H\hookrightarrow W_{n,m-1}$ and $H\hookrightarrow G$. Then
$W_{n,m-1}\to G_0$ and $G\to G_0$ are closed immersions, and $G_0/W_{n,m-1}\cong\alpha_p$. Thus
we need only to prove that $G_0$ can be embedded into $W_{n,m}^2$ as a closed subgroup scheme.

Note that $W_{n,m-1}$ is a closed subgroup scheme of $W_{n,m}$. Let Бо $G_1$ be the
push-out of $W_{n,m-1}\to W_{n,m}$ and $W_{n,m-1}\to G_0$. Then there is a commutative diagram
with exact rows:
$$\CD 0\to @. W_{n,m-1} @>>> G_0 @>>> \alpha_p @. \to 0 \cr
@. @VVV @VVV @VV{\text{id}}V @. \cr
0\to @. W_{n,m} @>>> G_1 @>{\lambda}>> \alpha_p @. \to 0 \endCD\tag 3.10$$
Let $G_i=\text{Spec}(A_i)$, $W_{n,m-1+i}=\text{Spec}(B_i)$ ($i=0,1$). Since $k$ is perfect,
as a $k$-algebra $A_1$ is of the form $k[y_1,...,y_r]/(y_1^{p^{n_1}},...,y_r^{p^{n_r}})$
by Fact 1.1.iii), and $B_1\cong k[x_0,...,x_n]/(x_0^{p^m},...,x_n^{p^m})$. Hence
there are only two possible cases:
either
$$A_1\cong k[y_0,...,y_n]/(y_0^{p^{m+2}},y_1^{p^{m+1}},...,y_n^{p^{m+1}}) \tag 3.11$$
or
$$A_1\cong k[y_0,...,y_{n+1}]/(y_0^{p^{m+1}},...,y_n^{p^{m+1}},y_{n+1}^p) \tag 3.12$$
If (3.11) holds,
it is easy to see that $\text{ker}(F^{m+1}_{G_1/k})$ has length $(m+1)(n+1)$, and
$W_{n,m}\subset\text{ker}(F^{m+1}_{G_1/k})$ also has length $(m+1)(n+1)$, hence
$W_{n,m}=\text{ker}(F^{m+1}_{G_1/k})$. Note that obviously $F_{G_0/k}^{m+1}=0$, hence
$G_0\subset\text{ker}(F^{m+1}_{G_1/k})=W_{n,m}$.

Taking the Cartier dual of (3.10) we get a commutative diagram with exact rows:
$$\CD 0\to @. \alpha_p @>>> G_1^D @>>> W_{n,m}^D @. \to 0 \cr
@. @VV{\text{id}}V @VV{q}V @VV{q'}V @. \cr
0\to @. \alpha_p @>>> G_0^D @>>> W_{n,m-1}^D @. \to 0 \endCD\tag 3.13$$
in which $q$ and $q'$ are faithfully flat. Note that
$\text{ker}(F_{W_{n,m}^D/k})\cong\text{coker}(V_{W_{n,m}/k})^D$ has length $m+1$
and $F_{W_{n,m}^D/k}^{n+1}=0$, while $W_{n,m}^D$ has length $(m+1)(n+1)$, hence
$B_1^D\cong k[y_0,...,y_m]/(y_0^{p^{n+1}},...,y_m^{p^{n+1}})$ as a $k$-algebra.
Thus the structure of $A_1^D$ has only two possible cases: either
$$A_1^D\cong k[y_0,...,y_m]/(y_0^{p^{n+2}},y_1^{p^{n+1}},...,y_m^{p^{n+1}}) \tag 3.14$$
or
$$A_1^D\cong k[y_0,...,y_{m+1}]/(y_0^{p^{n+1}},...,y_m^{p^{n+1}},y_{m+1}^p) \tag 3.15$$
If (3.14) holds,
it is easy to see that $B=\text{Spec}k[\bar y_0^p,\bar y_1,...,\bar y_m]\subset A_1^D$ gives $G_1^D$ a
quotient group scheme structure $\bar G\cong\text{ker}(V_{G_1/k}^{n+1})^D$, whose length is $(m+1)(n+1)$,
and obviously $B_1^D\subset B$, hence the projection $G_1^D\to\bar G$ induces $\bar G\cong W_{n,m}^D$.
Note that obviously $V_{G_0/k}^{n+1}=0$, hence $A_0^D\subset B$, therefore the projection $q$ factors
through $W_{n,m}^D$, showing that $G_0\to G_1$ factors through $W_{n,m}$.

The remaining case is that (3.12) and (3.15) hold.
Let $H_0=\text{ker}(F_{G_1/k})$. Then $H_0$ has length $n+2$, and $H_0^D$ satisfies
$F_{H_0^D/k}^{n+1}=0$. Let $R'$ be the structure ring of $H_0^D$. Then
$R'\cong k[y_0,y_1]/(y_0^{p^{n+1}},y_1^p)$. Since
$V_{H_0^D/k}=0$, by Fact 1.2.iii) we can take $y_0,y_1$ to be elements in $\alpha (H_0^D)$,
hence $k[\bar y^{p^n},\bar y_1]\subset R'$ gives a quotient group scheme $H_1$ of
$H_0^D$ (in fact $H_1=\text{coker}(F_{H_0^D/k})$, even without Fact 1.2 one can see
$\text{coker}(F_{H_0^D/k})$ has length 2, because $\text{ker}(F_{H_0^D/k})$ has length 2). Note that
$V_{H_1/k}=0$ and $F_{H_1/k}=0$, hence $H_1\cong\alpha_p\times_k\alpha_p$ by Fact 1.2.v).
Taking Cartier dual we get a closed subgroup scheme $H_1^D\cong\alpha_p\times_k\alpha_p$
of $H_0$. But by $a(W_{n,m})=1$ we see that $H_1^D\cap W_{n,m}\cong\alpha_p$, hence the projection
$H_1^D\to G_1\buildrel{\lambda}\over\to\alpha_p$ is surjective. This gives a section of $\lambda$,
thus $G_1\cong W_{n,m}\times_k\alpha_p$. \ \ \ Q.E.D.

\vskip 10pt
We now see several consequences of this proposition.

First, noting that the structure ring of $\Cal W_{n,m}$ is generated by
$x_0,...,x_n\in D(\Cal W_{n,m})$ as an $\Bbb F_p$-algebra, we see $R$ is generated by $D(G)$
as a $k$-algebra by Proposition 3.1.

Second, if $a(G)=1$, letting $H\subset G$ be the unique subgroup scheme isomorphic to $\alpha_p$, we
can take $x\in D(G)$ such that the image of $x$ in $D(H)$ is non-zero, thus $\text{ker}(f_x)$ does not
contain any subgroup scheme isomorphic to $\alpha_p$, hence $\text{ker}(f_x)=0$, i.e. $f_x$ is a
closed immersion; conversely, if there exists
$x\in D(G)$ such that $f_x$ is a closed immersion, then by Corollary 2.3.iii) we see that
$a(G)=1$. In this case $R$ is generated by $x$, $V^*x$, ..., $V^{n*}x$ as a $k$-algebra.

Therefore if $a(G)=1$ and $F_{G/k}^{m+1}=0$, $V_{G/k}^{n+1}=0$, then $G$ can be embedded into
$W_{n,m}$ as a closed subgroup scheme. In this case if $G$ has length $(m+1)(n+1)$, then the
embedding $G\to W_{n,m}$ is an isomorphism. In particular, take
$k=\Bbb F_p$, $G=\Cal W_{m,n}^D$, then $F_{G/k}^{m+1}=0$,
$V_{G/k}^{n+1}=0$ and $G$ has length $(m+1)(n+1)$, hence $\Cal W_{m,n}^D\cong\Cal W_{n,m}$.
In the next section we will specialize this isomorphism.

Third, we show that if $i:H\hookrightarrow G$ is a closed subgroup scheme, then the $A$-module
homomorphism $i^*:D(G)\to D(H)$ is surjective. For any
$\alpha_p\cong H'\subset H$, we can take $x\in D(G)$ such that the image $\bar x$ of $x$
in $D(H')$ is non-zero. Denote $x'=i^*(x)\in D(H)$. For any $y\in D(H)$, if its image in $D(H')$
is $\bar y=c\bar x$ ($c\in k$), then by Proposition 3.1.iii) we see that
$y\dotplus (-cx')\in\text{ker}(D(H)\to D(H'))= D(H/H')$. By induction on the length of
$G$, we may assume $D(G/H')\to D(H/H')$ is surjective, hence
there exists $z'\in D(G/H')\subset D(G)$ such that $i^*(z')=y\dotplus (-cx')$.
Let $z=z'\dotplus cx$, then we have
$$i^*(z)=i^*z'\dotplus ci^*x=y\dotplus (-cx')\dotplus cx'=y \tag 3.16$$
From this and Proposition 3.1.iii) we see that
$D$ is an exact functor. Furthermore, obviously $D$ preserves direct sums, and by definition it is
obvious that for any $G,G'\in\text{Ob}(\goth A\goth b_{\text{inf}}^{\text{inf}})$, the map
$Hom_k(G',G)\to Hom_A(D(G),D(G'))$ given by $D$ is a homomorphism of additive groups, hence
$D$ is an additive and exact contravariant functor of abelian categories.

In particular, if $k=\Bbb F_p$, noting that the powers of $V$ gives a filtration of $\Cal W_{n,m}$
whose factors are all isomorphic to
$\Bbb G_{a/k}[F^{m+1}]$, by Fact 1.2.iii) and inducition we see that $D(\Cal W_{n,m})$ is generated
by $x_n$ as an $A$-module. Thus it is easy to see that
$$D(\Cal W_{n,m})\cong A/A(F^{m+1},V^{n+1}) \tag 3.17$$
Hence for any $G\in\text{Ob}(\goth A\goth b_{\text{inf}}^{\text{inf}})$,
$D(G)$ is a finitely generated $A$-module such that the actions of $F$ and $V$ are both nilpotent,
and $a(G)=1$ if and only if $D(G)$ is generated by one element as an $A$-module.

By Fact 1.2.ii) we see $D(G)$ has a filtration whose factors are all isomorphic to
$D(\alpha_p)\cong k$, hence the length of $G$ is equal to $l_A(D(G))$. Summarizing we have

\vskip 10pt
\flushpar{\bf Corollary 3.1}. Let $k$ be a perfect field of characteristic $p>0$.

\vskip 5pt\leftskip=30pt\parindent=-10pt
i) If $G=\text{Spec}(R)\in\text{Ob}(\goth A\goth b_{\text{inf}}^{\text{inf}})$,
then $R$ is generated by $D(G)$ as a $k$-algebra, and $l_A(D(G))=l(G)$.

ii) $D:\goth A\goth b_{\text{inf}}^{\text{inf}}\to\goth M_A$ is an additive and exact
contravariant functor of abelian categories.

iii) For $G=\text{Spec}(R)\ne 0\in\text{Ob}(\goth A\goth b_{\text{inf}}^{\text{inf}})$, the following
conditions are equivalent:

\vskip 5pt\leftskip=50pt\parindent=-10pt
1) $a(G)=1$;

2) $D(G)$ is generated by one element as an $A$-module;

3) There exists $x\in D(G)$ such that $R$ is generated by all $V^{i*}x$'s as a $k$-algebra;

4) $G$ is isomorphic to a closed subgroup scheme of some $W_{n,m}$.

\vskip 5pt\leftskip=30pt\parindent=-10pt

iv) If $G\in\text{Ob}(\goth A\goth b_{\text{inf}}^{\text{inf}})$ santisfies $a(G)=1$ and
$F_{G/k}^{m+1}=0$, $V_{G/k}^{n+1}=0$, $l(G)=(m+1)(n+1)$, then $G\cong W_{n,m}$. In
particular $\Cal W_{m,n}^D\cong\Cal W_{n,m}$.

v) $D(W_{n,m})\cong A/A(F^{m+1},V^{n+1})$. This means that the elements of $D(W_{n,m})$
are the sums of some $c_{ij}x_i^{p^j}$ ($0\le i\le n$, $0\le j\le m$,
$c_{ij}\in k$) under the addition $\dotplus$.

\vskip 10pt\leftskip=0pt\parindent=20pt
By Corollary 3.1.iv) we see another basic property of $W_{n,m}$:
$$\dim_kLie(W_{n,m}/k)=n+1,\ \ \dim_kLie(W_{n,m}^D/k)=m+1 \tag 3.18$$

Let $\goth M_A^0\subset\goth M_A$ be the full subcategory of all finitely generated
$A$-modules such that the actions of $F$ and $V$ are both nilpotent. The above corollary
already shows that the restriction $D:\goth A\goth b_{\text{inf}}^{\text{inf}}\to\goth M_A^0$
is an anti-equivalence of abelian categories. But we prefer to give $D^{-1}$ directly
as follows.

For any $M\in\text{Ob}(\goth M_A^0)$, take $n$, $m$ such that $F^{m+1}=V^{m+1}=0$ on $M$,
thus get an $A$-module epimorphism $g:(A/A(F^{m+1},V^{n+1}))^{\oplus r}\to M$ for some $r$.
Let $R_0$ be the structure ring of $G_0=W_{n,m}^{\oplus r}$ and identify $D(G_0)$
with $(A/A(F^{m+1},V^{n+1}))^{\oplus r}$. By Definition 3.1 it is easy to see that
the ideal of $R_0$ generated by $\text{ker}(g)$ defines a closed subgroup scheme $G\subset G_0$.
By Corollary 3.1.ii) we have $D(G)\cong M$, and it is easy to check that $G$ is uniquely determined
by $M$ up to isomorphism. Hence we can denote $G=D^{-1}(M)$.

Let $f:M\to M'$ be a morphism in $\goth M_A^0$. Take an $A$-module epimorphism
$g:(A/A(F^{m+1},V^{n+1}))^{\oplus r}\to M$ as above, defining a closed immersion
$G=D^{-1}(M)\hookrightarrow G_0$. Let $G'=D^{-1}(M')$, the by Corollary 3.1.ii) we see that
$h=f\circ g:(A/A(F^{m+1},V^{n+1}))^{\oplus r}\to M'$ induces a group scheme homomorphism
$\eta :G'\to G_0$. By $h(\text{ker}(g))=0$ we see that $\eta$ induces a group scheme
homomorphism $\phi :G'\to G$. It is easy to check that $\phi$ is uniquely determined
by $f$ up to isomorphism, and $D(\phi )=f$. Hence we can denote $\phi =D^{-1}(f)$.

The above procedure gives the inverse $D^{-1}$ of $D$, showing
$D$ is an anti-equivalence of categories. Note that $D$ is an anti-equivalence of abelian
categories, because $D$ and $D^{-1}$ preserves direct sum (product), kernel and cokernel
by their universalities. Summarizing we have

\vskip 10pt
\flushpar{\bf Theorem 3.1}. For any perfect field $k$ of characteristic $p>0$,
let $\goth M_A^0\subset\goth M_A$ be the full subcategory of all finitely generated
$A$-modules such that the actions of $F$ and $V$ are both nilpotent. Then the
restriction of the Dieudonn\'e module functor
$D:\goth A\goth b_{\text{inf}}^{\text{inf}}\to\goth M_A^0$ is an anti-equivalence of abelian
categories, whose inverse $D^{-1}$ is given by the above constructive procedure.

\vskip 10pt
From the above arguments we also see that

\vskip 10pt
\flushpar{\bf Corollary 3.2}. Let $G=\text{Spec}(R)\in\text{Ob}(\goth A\goth b_{\text{inf}}^{\text{inf}})$.
Then an ideal $I\subset R$ defines a closed subgroup scheme if and only if $I$ is generated by
Dieudonn\'e elements (in it).

\vskip 30pt
\flushpar {\bf 4. Duality and quasi-polarizations}

\vskip 10pt
In \S 3 we see that $\Cal W_{m,n}^D\cong\Cal W_{n,m}$. In this section we give a
``standard'' isomorphism $\Cal W_{m,n}^D\buildrel{\simeq}\over\longrightarrow\Cal W_{n,m}$
explicitly. For this we need to give a Dieudonn\'e element in $A_{n,m}^D$, the structure
ring of $\Cal W_{m,n}^D$.
We go in an indirect way: first choose a Dieudonn\'e element of $\Cal W_{m,n}^D$, and then
determine it as an $\Bbb F_p$-linear functional on $A_{n,m}$.
For convenience denote $k=\Bbb F_p$ and $W_{n,m}=\Cal W_{n,m}$.
Again denote by $A_{n,m}$ the structure ring of $W_{n,m}$.

By Corollary 3.1, we know that there is a Dieudonn\'e element $y\in D(W_{n,m}^D)\subset A_{n,m}^D$
such that $A_{n,m}^D$ is generated by all $y$, $V^*y$,..., $V^{m*}y$ as a $k$-algebra.
Let $H\subset W_{n,m}^D$ be the closed subgroup scheme isomorphic to $\alpha_p$. Then
$H^D$ is canonically isomorphic to the quitient group scheme of $W_{n,m}$
isomorphic to $\alpha_p$, whose Dieudonn\'e generator (i.e. gnenrator of its $\alpha$-module) can be
chosen to be $V^{n*}x_n^{p^m}=x_0^{p^m}$. Hence (replacing $y$ by $cy$ for some $c\in\Bbb F_p^*$)
we can choose $y$ so that $y(x_0^{p^m})=1$.

Note that $W_{0,0}\cong\alpha_p$, and by Fact 1.2.i) we have
$D(\alpha_p^D)\cong Lie(\alpha_p/\Bbb F_p)$, whose generator $y$ satisfies
$y(x_0^i)=\delta_{i1}$ as a linear functional on $A_{0,0}$. Denote by $D_y\in Lie(\alpha_p/k)$
the left invariant derivation corresponding to $y$. Then by
$a^*(x_0)=x_0\otimes 1+1\otimes x_0$ it is easy to see that
$D_y(x_0^i)=ix_0^{i-1}$, i.e. $D_y=\frac{\text{\rm d}}{\text{\rm d}x_0}$.

In general case, note that
$$\aligned & V^{i*}y(x_0^{p^{m-i}})=y(F^{i*}V^{n*}x_n^{p^{m-i}})=y(V^{n*}x_n^{p^m})=1
\tpt (0\le i\le m) \cr 
& F^{i*}y(x_{n-i}^{p^m})=y(V^{i*}V^{(n-i)*}x_n^{p^m})=y(V^{n*}x_n^{p^m})=1
\tpt (0\le i\le n) \endaligned\tag 4.1$$
and $V^{i*}y$ (resp. $F^{i*}y$) can be viewed as an element of  $A_{n,m-i}^D$ (resp. $A_{n-i,m}^D$).
If $y(V^{n*}x_n^{p^{m-1}})=c$, replacing $y$ by $y\dotplus (-cV^*y)$ we can reduce to the
case $y(V^{n*}x_n^{p^{m-1}})=0$; similarlay, if $y(V^{(n-1)*}x_n^{p^m})=c$, replacing $y$ by
$y\dotplus (-cF^*y)$ we can reduce to the case $y(V^{(n-1)*}x_n^{p^m})=0$. By induction we
can eventually reduce to the case
$$y(x_{n-i}^{p^j})=y(V^{i*}x_n^{p^j})=(V^{r*}y)(x_{n-i}^{p^{j-r}})=\delta_{in}\delta_{jm}
\tpt (\forall i\le n,j,r\le j) \tag 4.2$$
In the following we will show that, in this case
$$y(\prod_{i=0}^nx_i^{j_i})=\cases 1 & \text{ if } j_1=...=j_n=0,\ j_0=p^m \cr
0 & \text{ otherwise} \endcases\tag 4.3$$
as a $k$-linear functional.

For a monomial
$$\alpha =x_0^{i_0}\cdots x_n^{i_n}\in A_{n,m} \tag 4.4$$
define its ``degree'' to be
$$d(\alpha )=i_0+pi_1+\cdots+p^ni_n \tag 4.5$$
For any $d\in\Bbb N$, let $A_{n,m}^d\subset A_{n,m}$ be the
$k$-linear subspace generated by all monomials of degree $d$. It is easy to see this defines
a graded $k$-algebra structure on $A_{n,m}$. By Corollary 2.1.iv) we see that
every monomial $\alpha\otimes\beta$ in
$a^*(x_n)=\phi_n(x_0\otimes 1,1\otimes x_0,...,x_n\otimes 1,1\otimes x_n)$ satisfies
$d(\alpha )+d(\beta )=p^n$, hence $a^*$ is a graded homomorphism.

For a monomial $\alpha$ in (4.4), we can denote
$$a_{W_{n,m}}^*(\alpha )=\sum\limits_{0\le i_0,...,i_n<p^{m+1}}x_0^{i_0}\cdots x_n^{i_n}\otimes\alpha_{i_0,...,i_n}
\tag 4.6$$
where $\alpha_{i_0,...,i_n}\in A_{n,m}$. Suppose (4.3) holds. Then by (4.6) we have
$$D_y(\alpha )=\alpha_{p^m,0,...,0} \tag 4.7$$
From this we see that every monomial of $D_y(\alpha )$ has degree equal to
$\deg (\alpha )-p^m$. Furthermore, by (4.1) and (4.3) we can get
$$V^{r*}y(\prod_{i=0}^nx_i^{j_i})=\cases 1 & \text{ if } j_1=...=j_n=0,\ j_0=p^{m-r} \cr
0 & \text{ otherwise} \endcases\tag 4.8$$
From this we can see that every term of $D_{V^{r*}y}(\alpha )$ has degree equal to
$\deg (\alpha )-p^{m-r}$. Therefore we have
$$D_{V^{r*}y}(A_{n,m}^d)\subset A_{n,m}^{d-p^{m-r}} \tpt (\forall d,r) \tag 4.9$$
In particular, when $d<p^{m-r}$ we have $D_{V^{r*}y}(A_{n,m}^d)=0$.

We now show (4.3). First note that if $i:H=\text{Spec}(R_1)\subset G=\text{Spec}(R)$ is a closed
immersion in $\goth A\goth b_{\text{inf}}^{\text{inf}}$ and $I=\text{ker}(i^*)\subset R$, then by definition
we know that $i^{D*}=i^{*D}:R_1^D\to R^D$ is a $k$-algebra monomorphism, whose image is
$$\{\beta\circ i^*|\beta\in R_1^D\}=\{\alpha\in R^D|\alpha (I)=0\} \tag 4.10$$
In particular, let $G=W_{n,m}$ and $H\subset G$  be the image of $F_{G/k}$,
then $H=W_{n,m-1}$ and the Cartier dual of $i:H\to G$ can be viewed as the coimage
of $V_{G^D/k}$, which is isomorphic to the image of $i^D:G^D\to H^D$. Hence
$$V^*(A^D)=\{\alpha\in R^D|\alpha ((x_0^{p^m},...,x_n^{p^m}))=0\} \tag 4.11$$
By the definition of Dieudonn\'e element we have
$$a_{W_{n,m}^D}^*(y)=\phi_m(V^{m*}y\otimes 1,1\otimes V^{m*}y,...,y\otimes 1,1\otimes y)
\tag 4.12$$
Note that $\mu =a_{W_{n,m}^D}^{*\vee}:A_{n,m}\otimes A_{n,m}\to A_{n,m}$ (the dual of
$a_{W_{n,m}^D}^*$) is the multiplication of $A_{n,m}$, and the left hand side of (4.12) is
equal to $y\circ\mu$. Hence for any $x,x'\in A_{n,m}$ we have
$$D_y(xx')=\mu\circ\phi_m(D_{V^{m*}y}\otimes 1,1\otimes D_{V^{m*}y},...,D_y\otimes 1,1\otimes D_y)
(x\otimes x') \tag 4.13$$
by Fact 1.1.i). Let $x=x_0^{p^m}$. Then by (4.2) and (4.13) we can see that
$$D_y(x_0^{p^m}x')=x'+x_0^{p^m}D_y(x') \tag 4.14$$
Similarly, let $x=x_i^{p^m}$ ($i>0$). Then we get
$$D_y(x_i^{p^m}x')=x_i^{p^m}D_y(x') \tpt (\forall i>0)  \tag 4.15$$
By (4.14) and (4.15), using induction we get
$$D_y((\prod_{i=0}^nx_i^{j_i})^{p^m}x')
=\cases j_0(x_0^{j_0-1}\prod_{i=1}^nx_i^{j_i})^{p^m}x'+(\prod_{i=0}^nx_i^{j_i})^{p^m}D_y(x')
& \text{ if } j_0>0 \cr
(\prod_{i=0}^nx_i^{j_i})^{p^m}D_y(x') & \text{ if } j_0=0 \endcases\tag 4.16$$
In particular, for any monomial $x'=\prod\limits_{i=0}^nx_i^{j_i}$ and any $j>0$ we have
$y(x_j^{p^m}x')=\overline{D_y(x_j^{p^m}x')}=0$, and when $x'\ne 1$ we have
$y(x_0^{p^m}x')=\overline{D_y(x_0^{p^m}x')}=0$. Thus for (4.3), the remaining case
is that $j_0,...,j_n<p^m$.

Next we use induction on $m$. Let $y'=V^*y$, $\alpha =\prod\limits_{i=0}^nx_i^{j_i}$ ($j_0,...,j_n<p^m$).
By induction hypothesis we have
$$V^*y(\alpha )=\cases 1 & \text{ if } j_1=...=j_n=0,\ j_0=p^{m-1} \cr
0 & \text{ otherwise} \endcases\tag 4.17$$
Hence
$$y(\alpha^p)=\cases 1 & \text{ if } j_1=...=j_n=0,\ j_0=p^{m-1} \cr
0 & \text{ otherwise} \endcases\tag 4.18$$
The remaining case is that $\alpha$ is of the form $\prod\limits_{i=0}^nx_i^{j_i}$ ($j_0,...,j_n<p^m$)
and at least one $j_i$ is not divisible by $p$. Let $x=x_i$, $x'=x_i^{j_i-1}\prod\limits_{l\ne i}x_l^{j_l}$, and we may assume $x'\ne 1$. It is easy to see that the right hand side of (4.13) is not in
$(x_0,...,x_n)$ if and only if the right hand side of (4.12) has a term
$D_{V^{i*}y}^{p^n}\otimes D'$ such that $D'(x')\in k-\{ 0\}$. Such a term has degree
$\ge p^{n+m-i}$. But by Corollary 2.1.iv) we know that such a term has degree $p^m$. If $n\ge m$, this can
happen only when $n=m$, $D'=1$ and $x'\in k$. Hence the right hand side of (4.13) is always in
$(x_0,...,x_n)$ when $n\ge m$, thus $y(xx')=\overline{D_y(xx')}=0$. 

We now use induction on $n$. Let $y'=y^p$. By (4.2) we have
$$y'(x_{n-i}^{p^j})=\delta_{i(n-1)}\delta_{jm} \tpt (0\le i\le n-1,\ 0\le j\le m) \tag 4.19$$
Note that the quotient group scheme of $W_{n,m}^D$ defined by $y'\in A_{n,m}^D$ is the image
of $F_{G^D/k}$, whose Cartier dual is the closed subgroup scheme of $W_{n,m}$ given by the
ideal $(x_0)$, hence whose structure ring is generated by $x_1,...,x_n$. By induction
hypothesis we have
$$y'(\prod_{i=1}^nx_i^{j_i})=\cases 1 & \text{ if } j_2=...=j_n=0,\ j_1=p^m \cr
0 & \text{ otherwise} \endcases\tag 4.20$$
Note that $y'(\alpha )=y(V^*\alpha )$, hence (4.20) gives
$$y(\prod_{i=0}^{n-1}x_i^{j_i})=\cases 1 & \text{ if } j_1=...=j_{n-1}=0,\ j_0=p^m \cr
0 & \text{ otherwise} \endcases\tag 4.21$$
The remaining case is that $j_n>0$. Let $x=x_n$. Then for any monomial $x'$, by above we see that
the right hand side of (4.13) is not in $(x_0,...,x_n)$ if and only if the right hand side of (4.12)
has a term $D_{V^{m*}y}^{p^n}\otimes D'$ such that $D'(x')\in k-\{ 0\}$. Such a term has degree
$\ge p^n$ and $=p^m$, if $n\ge m$ this can only happen when $n=m$, $D'=1$
and $x'\in k$. Hence when $n\ge m$ the right hand side of (4.13) is always in $(x_0,...,x_n)$, thus
$y(xx')=\overline{D_y(xx')}=0$.

Therefor we can show that (4.3) holds whenever $m\le n$. But for general $n,m$, since
$W_{n,m}$ can be embedded into $W_{n+m,m}$, it ie easy to see (4.3) still holds.
Summarizing we have

\vskip 10pt
\flushpar{\bf Theorem 4.1}. Denote by $A_{n,m}$ the structure ring of $\Cal W_{n,m}$.
Then the linear functional $y\in A_{n,m}^D$ defined by (4.3) is a Dieudonn\'e generator
of $\Cal W_{n,m}^D$. Furthermore,

\vskip 5pt\leftskip=30pt\parindent=-10pt
i) For any $0\le i\le n$ and $0\le j\le m$, $z=V^{j*}y^{p^i}\in D(\Cal W_{n,m}^D)$ is the only
Dieudonn\'e element of $\Cal W_{n,m}^D$ satisfying
$$z(x_r^{p^{m-s}})=\delta_{ir}\delta_{js} \tpt (\forall r\le n,\ j\le m) \tag 4.22$$

ii) (4.14) and (4.15) hold for any $x'\in A_{n,m}$.

iii) $A_{n,m}$ has a graded $\Bbb F_p$-algebra structure
$$A_{n,m}=\bigoplus\limits_{d=0}^{(p^{n+1}-1)(p^{m+1}-1)/(p-1)}A_{n,m}^d \tag 4.23$$
in which every monomial $\alpha =x_0^{i_0}\cdots x_n^{i_n}$ is a homogeneous element of degree
$d(\alpha )=i_0+pi_1+\cdots+p^ni_n$, and $A_{n,m}^d$ is the $\Bbb F_p$-linear subspace generated
by monomials of degree $d$. For any $0\le i\le n$ and $0\le j\le m$ we have
$$D_{V^{j*}y}^{p^i}(A_{n,m}^d)\subset A_{n,m}^{d-p^{i+m-j}} \tpt (\forall d) \tag 4.24$$
In particular $D_{V^{j*}y}^{p^i}(A_{n,m}^d)=0$ when $d<p^{i+m-j}$.

\vskip 10pt\leftskip=0pt\parindent=20pt
Let $k$ be a perfect field of characteristic $p>0$. Still denote $W_{n,m}=\Cal W_{n,m}\otimes k$,
$W=\Cal W\otimes k$, and $A$ the ring in (3.8).
For any $G\in\text{Ob}(\goth A\goth b_{k\text{inf}}^{\text{inf}})$,
we now study the relationship of $D(G)$ and $D(G^D)$ as $A$-modules.

Take any monomorphism $i:G^D\hookrightarrow W_{n,m}^r$. Then
$i^D:W_{n,m}^{Dr}\to G$ is an epimorphism, and any homomorphism $G^D\to W$ factors
through $W_{n,m}$. Fix the standard isomorphism $W_{n,m}^D\cong W_{m,n}$ by Theorem 4.1.
Then an element $y\in D(G^D)$ corresponds to a homomorphism
$y^D: W_{m,n}\to G$, and this corresponds to an $A$-module homomorphism
$y^{D*}:D(G)\to D(W_{m,n})\cong A/(F^{n+1},V^{m+1})$ by Corollary 3.1.
This gives a mapping $f:D(G^D)\to Hom_A(D(G),A/(F^{n+1},V^{m+1}))$.
By the surjectivity of $i^D$ we see $f$ is injective; By Corollary 3.1 we see
that any element of $Hom_A(D(G),A/(F^{n+1},V^{m+1}))$ corresponds to a
homomorphism $W_{m,n}\to G$, hence $f$ is surjective also. It is easy to see that for
any $a\in k$ and $x\in D(G)$, we have $\tau (a)\dot\times x=ax$, hence
$f(\tau (a)\dot\times x)=\tau (a)f(x)$. Furthermore, it is easy to see that
$f$ is additive, and
$$f(Fy)(x)=(f(y)(Vx))^{\sigma},\ f(Vy)(x)=(f(y)(Fx))^{\sigma^{-1}} \tpt (\forall x\in D(G),\ 
y\in D(G^D)) \tag 4.25$$
Thus $f(py)=pf(y)$, therefore $f$ is a $W(k)$-linear isomorphism. Let $r=\min (m,n)$.
For any $x\in D(G)$ and $y\in D(G^D)$, let $\langle x,y\rangle\in W_r(k)$ be the
constant term of $f(y)(x)\in A/(F^{n+1},V^{m+1})$. In this way we express $f$ as a
$W(k)$-bilinear form $\langle ,\rangle :D(G)\times D(G^D)\to W_r(k)$
which satisfies
$$\langle Fx,y\rangle =\langle x,Vy\rangle^{\sigma} ,
\ \ \langle Vx,y\rangle =\langle x,Fy\rangle^{\sigma^{-1}} \tpt (\forall x\in D(G),y\in D(G^D))
\tag 4.26$$
It is a perfect paring over $W(k)$.

We omit more facts about duality because they are more or less well-known. The only difference from the
known formulation of duality is that $D(G)$ is a subset of the structure ring of $G$ here.

\vskip 30pt
\flushpar {\bf 5. The connection with differential operators}

\vskip 10pt
On the polynomial ring $\Bbb Z[x]$ we can define differential operators
$D^{(n)}=\frac{1}{n!}\frac{\text{\rm d}^n}{\text{\rm d}x^n}$ for any $n\in\Bbb N$.
Indeed, for any $\phi\in\Bbb Z[x]$, $D^{(n)}(\phi )$ is the coefficient of
$\text{\rm d}x^n$ in the expansion of $\phi (x+\text{\rm d}x)$ as a polynomialof
$\text{\rm d}x$.
Therefore such differential operators $D^{(n)}$ can also be defined on $\Bbb F_p[x]$.

Note that since $x_0,...,x_n$ are independent variables over $\Bbb F_p$ in the structure ring of
$\Cal W_n$, if we denote the differential operators
$$D_{i,j}=\frac{1}{p^{m-j}!}\big(\frac{\partial}{\partial x_i}\big)^{p^{m-j}}\big|_0 \tag 5.1$$
then Theorem 4.1.i) says that the $\Bbb F_p$-linear functional $V^{j*}y^{p^i}$ on $A_{n,m}$ can be
understood as $D_{i,j}$ ($\forall 0\le i\le n$, $0\le j\le m$).

\vskip 10pt
\flushpar{\bf Example 5.1}. In the case $n=1$, by direct calculation we can get
$$\phi_1(x_0,y_0,x_1,y_1)=x_1+y_1-\sum_{i=1}^{p-1}\frac{(p-1)!}{i!(p-i)!}x_0^iy_0^{p-i} \tag 5.1$$
Since we are studying $\Bbb F_p$-schemes, we often use $(p-1)!\equiv -1$ (mod $p$) to rewrite
the right hand side of (5.1) as
$$x_1+y_1+\sum_{i=1}^{p-1}\frac 1{i!(p-i)!}x_0^iy_0^{p-i} \tag 5.2$$
For any $m$, by Theorem 4.1 we see that $D(\Cal W_{1,m}^D)$ has generators $y_j=V^{(n-j)*}y$
($0\le j\le m$), which can be understood as $D_{0,j}$. If $m=0$, by (5.2) it is not hard to
get
$$D_{y_0}=\frac{\partial}{\partial x_0}+x_0^{p-1}\frac{\partial}{\partial x_1},\ 
D_{y_0}^p=\frac{\partial}{\partial x_1} \tag 5.3$$

\vskip 10pt
For more details about the relations of Dieudonn\'e modules with differential operators,
we need some preparations. Note that the above differential operators
$D^{(n)}$ ($\forall n$) are commutative to each other.

\vskip 10pt
\flushpar{\bf Lemma 5.1}. Let $R=\Bbb Z[x]$ and $p$ be a prime number. Let $r\in\Bbb N$ and
$$D^{(p^i)}=\frac{1}{p^i!}\frac{\text{\rm d}^{p^i}}{\text{\rm d}x^{p^i}} \tpt (0\le i\le r) \tag 5.4$$
and denote $D=D^(1)$. Then

\vskip 5pt\leftskip=30pt\parindent=-10pt
i) For the ideal $I=(p,x^{p^r})\subset R$ and any $i<r$ we have $D^{(p^i)}(I)\subset I$, hence
$D^{(p^i)}$ induces a differential operator on $R/I$ (still denoted by $D^{(p^i)}$).
Furthermore,
$$\aligned \text{Diff}(R/I,\Bbb F_p) & =Hom(R/I,\Bbb F_p) \cr
& =\{ h\circ f(D^{(p^0)},D^{(p^1)},...,D^{(p^{r-1})})|
f\in\Bbb F_p[x_0,...,x_{r-1}]\} \endaligned\tag 5.5$$
where $h:R/I\to\Bbb F_p$ is the projection (modulo $(x)$), and we can take $f$ to be polonomials
whose exponents are all $<p$ (because $D^{(p^i)p}=0$ for each $i$).

ii) There exists a polynomial with $2(r+1)$ variables $\lambda_r\in\Bbb Z_{(p)}[x_0,y_0,...,x_r,y_r]$,
whose exponents are all $<p$, such that for any $a,b\in R$,
$$D^{(p^r)}(ab)=\lambda_r(D^{(p^0)}\otimes 1,
1\otimes D^{(p^0)},...,D^{(p^r)}\otimes 1,1\otimes D^{(p^r)})(a\otimes b) \tag 5.6$$

iii) For any $i\ge 0$ we have $D^{(p^r)}(\Bbb Z[x^{p^i}])\subset\Bbb Z[x^{p^i}]$, and the
restriction of $D^{(p^r)}$ on $\Bbb Z[x^{p^r}]$ modulo $p$ is a derivation.

\vskip 10pt\leftskip=0pt\parindent=20pt
\flushpar {\sl Proof}. We need to use several combinatorial facts. First, if
$m=i_0+i_1p+i_2p^2+\cdots +i_rp^r$ ($0\le i_0,...,i_r<p$), then
$$\frac{m!}{(p!)^{i_1}\cdots (p^r!)^{i_r}}\equiv i_0!i_1!\cdots i_r!\not\equiv 0
\ \ \ (\text{mod } p) \tag 5.7$$
Note that the left hand side of (5.7) is the coefficient of the monomial
$$\alpha =x_1\cdots x_{i_0}x_{i_0+1}^p\cdots x_{i_0+i_1}^p\cdots x_{i_0+\cdots +i_r}^{p^r}
\tag 5.8$$
in the polynimial $\phi =(x_1+\cdots +x_{i_0+\cdots +i_r})^m$, and
$$\phi\equiv (x_1+\cdots +x_{i_0+\cdots +i_r})^{i_0}(x_1^p+\cdots +x_{i_0+\cdots +i_r}^p)^{i_1}
\cdots (x_1^{p^r}+\cdots +x_{i_0+\cdots +i_r}^{p^r})^{i_r} \ \ \ (\text{mod } p) \tag 5.9$$
The coefficient of $\alpha$ in the right hand side is $i_0!i_1!\cdots i_r!$. Apply (5.7) to
$m=p^r-1=(p-1)(1+p+\cdots +p^{r-1})$ we get
$$1=\frac{p^r!}{(p^r-1)!p^r}\equiv (-1)^r\frac{p^r!}{(p!\cdots p^{r-1}!)^{p-1}p^r}
\ \ (\text{mod }p) \tag 5.10$$
Then by induction we get
$$\frac{p^r!}{p^{(p^r-1)/(p-1)}}\equiv (-1)^r\ \ (\text{mod }p) \tag 5.11$$
Applying (5.7) again, we see that for any $i<r$,
$$\frac{p^r!}{p^i!(p^r-p^i)!p^{r-i}}\equiv 1\ \ (\text{mod }p) \tag 5.12$$
By (5.11) we get
$$\frac{p^r!}{p^{r-1}!(p!)^{p^{r-1}}}\equiv 1\ \ (\text{mod }p) \tag 5.13$$
Then by induction we get
$$\frac{p^{r+s}!}{p^r!(p^s!)^{p^r}}\equiv 1\ \ (\text{mod }p) \tag 5.14$$

i) We need to check that $D^{(p^i)}(x^{p^r+m})\in I$ for any $m\ge 0$. It is enough
to consider the case when $m<p^i$.
We have $D^{(p^i)}(x^{p^r+m})=\frac{(p^r+m)!}{p^i!(p^r+m-p^i)!}x^{p^r+m-p^i}$. By (5.7) and (5.12) we have
$$\frac{(p^r+m)!}{p^i!(p^r+m-p^i)!}\equiv\frac{p^r!m!}{p^i!(p^r-p^i)!m!}\equiv 0
\ \ \ (\text{mod } p) \tag 5.15$$
Note also that for any $m=i_0+i_1p+i_2p^2+\cdots +i_{r-1}p^{r-1}$
($0\le i_0,...,i_{r-1}<p$), 
$$D^{(p^0)i_0}\circ\cdots\circ D^{(p^{r-1})i_{r-1}}(x^m)=
\frac{m!}{(p!)^{i_1}\cdots (p^{r-1}!)^{i_{r-1}}}\not\equiv 0 \ \ \ (\text{mod } p) \tag 5.16$$
and for any $j\ne m$ we have $D^{(p^0)i_0}\circ\cdots\circ D^{(p^{r-1})i_{r-1}}(x^j)\equiv 0$ (mod $(x)$).
Hence (5.5) holds.

ii) We have $D^{(p^r)}(ab)=\sum\limits_{m=0}^{p^r}\frac 1{m!(p^r-m)!}D^m(a)D^{p^r-m}(b)$,
hence it is enough to express $\frac 1{m!}D^m$ ($0<m<p^r$) as a polynimial of
$D^{(p^0)},...,D^{(p^{r-1})}$ with coefficients in $\Bbb Z_{(p)}$. Let
$m=i_0+i_1p+i_2p^2+\cdots +i_{r-1}p^{r-1}$ ($0\le i_0,...,i_{r-1}<p$). By (5.7) we see there
exists $c\in\Bbb Z_{(p)}^*$ such that
$$\frac 1{m!}D^m=c\frac 1{(p!)^{i_1}\cdots (p^{r-1}!)^{i_{r-1}}}D^m
=cD^{(p^0)i_0}\circ\cdots\circ D^{(p^{r-1})i_{r-1}} \tag 5.17$$
By (5.17) we can see the uniqueness of $\lambda_r$.

iii) The first assertion is obvious. For any multiple $m$ of $p^r$, decompose to $m=p^sq$,
where $s\ge r$ and $p\nmid q$. By (5.7) we get $\frac{m!}{(m-p^s)!p^s!}\equiv q$ (mod $p$),
hence by (5.12) we see that $\frac{m!}{(m-p^r)!p^r!}\equiv 0$ (mod $p$) when $s>r$. Therefore
$D^{(p^r)}(x^m)=\frac{m!}{p^r!(m-p^r)!}x^{m-p^r}$ is congruent to $\frac{m}{p^r}x^{m-p^r}$
modulo $p$ in any case. This shows that the restriction of $D^{(p^r)}$ on
$\Bbb F_p[x^{p^r}]\subset\Bbb F_p[x]$ is equal to
$\frac{\text{\rm d}}{\text{\rm d}y}$, where $y=x^{p^r}$. \ \ \ Q.E.D.

\vskip 10pt
View $\text{Spec}(R)$ as $\Bbb G_{a/\Bbb Z}$, whose addition is given by
$a^*(x)=x\otimes 1+1\otimes x$. For any $r\in\Bbb N$,
$\Bbb G_{a/\Bbb Z}\otimes\Bbb F_p\cong\Bbb G_{a/\Bbb F_p}$ has a closed subgroup scheme
$$H_r=\Bbb G_{a/\Bbb F_p}[F^{r+1}]\cong\alpha_{p^{r+1}}=\text{Spec}(B_r) \tag 5.18$$
where $B_r=\Bbb F_p[x]/(x^{p^{r+1}})$. Each $D^{(p^i)}$ ($0\le i\le r$) gives a left invariant
differential operator $t_i$ satisfying
$t_i(x^{p^j})=\delta_{ij}$ ($0\le i,j\le r$). Note that $x,x^p,...,x^{p^r}$ form a basis of the
$\alpha$-module of $H_r$. By Theorem 4.1 we see each $\bar t_i=t_i|_0\in B_r^D$ 
is a Dieudonn\'e element of $H_r^D$, where
$\bar t_r$ is a Dieudonn\'e generator, and $\bar t_i=V^{(r-i)*}\bar t_r$ ($0\le i\le r$).
By $V_{H_r/\Bbb F_p}=0$ we have $F_{H_r^D/\Bbb F_p}=0$, hence $t_i^p=0$ ($0\le i\le r$).
Therefore $B_r^D\cong\Bbb F_p[y_0,...,y_r]/(y_0^p,...,y_r^p)$, where each $y_i$ corresponds
to $\bar t_i$, and (5.6) gives
$$a^*(y_i)=\lambda_i(V^{i*}y_i\otimes 1,1\otimes V^{i*}y_i,V^{i-1*}y_i\otimes 1,
1\otimes V^{i-1*}y_i,...,y_i\otimes 1,1\otimes y_i)\ \ (0\le i\le r) \tag 5.19$$
This shows that $\lambda_i$ and $\phi_i$ take same values on $B_r^D$. Summarizing we have

\vskip 10pt
\flushpar{\bf Corollary 5.1}. The structure ring of $\alpha_{p^{r+1}}^D$ is isomorphic to
$\Bbb F_p[y_0,...,y_r]/(y_0^p,...,y_r^p)$, whose co-addition is given by (5.19).
Furthermore,
$$\lambda_i\equiv\phi_i\ \text{ (mod }(p,y_0^p,...,y_r^p)) \tpt (0\le i\le r) \tag 5.20$$
where $\phi_i$ is the polynomial defined in \S 2.

\vskip 10pt

\vskip 10pt
\flushpar{\bf Example 5.2}. In the case $r=1$, we have $\phi_1$ as in (5.1), while
$$\lambda_1(x_0,y_0,x_1,y_1)=x_1+y_1+\sum_{i=1}^{p-1}\frac{1}{i!(p-i)!}x_0^iy_0^{p-i} \tag 5.21$$
which is not in $\Bbb Z[x_0,y_0,x_1,y_1]$, but congruent to $\phi_1$ modulo $p$.

Let $k\supset\Bbb F_p$ be an algebraically closed field. It is easy to classify all
$G=\text{Spec}(A)\in\text{Ob}(\goth A\goth b_{\text{inf}}^{\text{inf}})$ of length 2
by Dieudonn\'e module theory. Up to isomorphism $D(G)$ has only four possible cases:
$$(A/A(F,V))^{\oplus 2},\ A/A(F^2,V),\ A/A(F,V^2),\ A/A(F-V,p) \tag 5.22$$
The first one corresponds to $\alpha_p^2$; the secong one corresponds to $\alpha_{p^2}$;
and the third one corresponds to $\alpha_{p^2}^D$, whose group scheme structure is explained in
Corollary 5.1. The last one corresponds to $E[p]$ for a supersingular elliptic curve
over $k$, whose group scheme structure can be calculated using Theorem 3.1.
We have $E[p]\cong\text{Spec}k[x_1]/(x_1^{p^2})$, where $x_1$ satisfies
$$a^*(x_1)=x_1\otimes_k1+1\otimes_kx_1+\sum_{j=1}^{p-1}\frac 1{j!(p-j)!}x_1^{jp}\otimes_kx_1^{(p-j)p},
\ 0^*(x_1)=0,\ (-1)^*(x_1)=-x_1 \tag 5.23$$

\vskip 30pt\parindent=0pt
\flushpar {\bf References}

\vskip 10pt
[BBM] P. Berthelot, L. Breen \& W.Messing: {\sl Th\'eorie de Dieudonn\'e cristalline
I\! I}. LNM 930, Springer-Verlag (1982).

\vskip 5pt
[BO] P. Berthelot \&\ A. Ogus: {\sl Notes on Crystalline Cohomology}. Princeton Univ.
Press \&\ Univ. Tokyo Press., Princeton (1978).

\vskip 5pt
[dJ1] J. de Jong: Moduli of abelian varieties and Dieudonn\'e modules of
finite group schemes (Thesis). Univ. of Utrecht (1992)

\vskip 5pt
[dJ2] A.J. de Jong: The moduli spaces of polarized abelian varieties.
Math. Ann. 295 (1993), 485-503.

\vskip 5pt
[Dem] M. Demazure: {\sl Lectures on p-divisible Groups}. LNM 302,
Springer-Verlag (1972).

\vskip 5pt
[DG] M. Demazure \&\ P. Gabriel - {\sl Groupes alg\'ebriques, I}. Masson, Paris
and North-Holland Pub. Co., Amsterdam (1970).

\vskip 5pt
[Ka] N. Katz: Slope filtration of F-crystals. Journ\'ees G\'eom.
Alg\'ebr.
Rennes 1978, Vol. I. Ast\'erisque 63, Soc. Math. France (1979), 113-164.

\vskip 5pt
[Laz] M. Lazard: {\sl Commutative Formal Groups}. LNM 443, Springer-Verlag,
(1975).

\vskip 5pt
[Li1] K. Li: Classification of supersingular abelian varieties. Math. Ann.
283 (1989), 333-351

\vskip 5pt
[Li2] Ke-Zheng Li: Differential Operators and automorphism schemes. Science China
Mathematics Volume 53 Number 9 (2010), 2363-2380

\vskip 5pt
[LO] Ke-Zheng Li \&\ Frans Oort: {\sl Moduli of Supersingular Abelian
Varieties}, LNM 1680. Springer (1998)

\vskip 5pt
[Ma] Yu. I. Manin: The theory of commutative formal groups over fields of
finite characteristic. Usp. Math. 18 (1963), 3-90; Russ. Math. Surveys 18
(1963), 1-80.

\vskip 5pt
[Me] W. Messing - {\sl The Crystals Associated to Barsotti-Tate Groups: with
Applications to Abelian Schemes}. LNM 264, Springer-Verlag (1972).

\vskip 5pt
[Mu1] D. Mumford: {\sl Lectures on Curves on an Algebraic Surface},
Annals of Math. Studies 59, Princeton University Press, Princeton (1966)

\vskip 5pt
[Mu2] D. Mumford: {\sl Abelian varieties}, 2nd print, Tata Inst. Fund.
Res. \& Oxford Univ. Press  (1974)

\vskip 5pt
[Oo1] F. Oort: {\sl Commutative group schemes}, LNM 15, Springer-Verlag (1966)

\vskip 5pt
[P] Richard Pink: Finite group schemes, Lecture course in WS 2004/05, ETH, Z\"urich

\vskip 5pt
[Ra] M. Raynaud: Sch\'emas en groupes de type $(p,..., p)$, Bull. Soc. math.
France 102 (1974), 241-280

\vskip 5pt
[S3.I] M. Demazure, A. Grothendieck et al - Schemas en Groups (SGA 3) I,
LNM 151, Springer-Verlag (1970)

\vskip 5pt
[Si] J. H. Silverman: {\sl The Arithmetic of Elliptic Curves}. GTM 106,
Springer-Verlag (1986).

\vskip 5pt
[SGA] A. Grothendieck et al.: {\sl S\'eminaire de G\'eom\'etrie Alg\'ebrique},
LNM225, 269, 270, 305. Springer-Verlag

\vskip 5pt
[TO] J. Tate \&\ F. Oort: Group schemes of prime order, Ann. Sc. Ecole Norm.
Sup. 3 (1970), 1-21

\vskip 20pt\parindent=0pt
Kezheng Li

Department of Mathematics

Capital Normal University

Beijing 100048, China

e-mail: kzli\@\!\, gucas.ac.cn

\bye